\documentclass[10pt,twocolumn,twoside]{IEEEtran}

\usepackage{amsmath,amssymb,amsthm}
\usepackage{graphicx}
\usepackage{xcolor}
\newcommand{\N}{\ensuremath{\mathbb{N}}}    
\newcommand{\R}{\ensuremath{\mathbb{R}}}    

\newcommand{\x}{\ensuremath{\mathbf{x}}}    
\newcommand{\y}{\ensuremath{\mathbf{y}}}    
\newcommand{\w}{\ensuremath{{\boldsymbol \omega}}}  

\newcommand{\rxyw}{\ensuremath{\R[\x,\y,\w]}}   
\newcommand{\rxy}{\ensuremath{\R[\x,\y]}}       
\newcommand{\rxw}{\ensuremath{\R[\x,\w]}}       
\newcommand{\rx}{\ensuremath{\R[\x]}}           

\renewcommand{\l}{\ensuremath{\ell}}            
\renewcommand{\d}{\ensuremath{\mathrm{d}}}      

\newcommand{\prob}{\ensuremath{\mathbb{P}_\w}}     
\newcommand{\st}{\ensuremath{\mbox{s.t.}}}      

\newcommand{\meas}{\ensuremath{M_+}}			
\newcommand{\vol}{\ensuremath{\mathrm{vol}}}	    

\newcommand{\dmu}{\mu(\d\w)}

\newcommand{\kxp}{\ensuremath{\mathcal{L}^{\x}}}
\newcommand{\kxpp}{\ensuremath{\mathcal{K}^{\x}}}
\newcommand{\kx}{\ensuremath{K^\x}}

\newcommand{\kxyw}{\ensuremath{K}}

\newcommand{\Y}{\ensuremath{\mathbf{Y}}}

\newtheorem{thm}{Theorem} 

\newtheorem{prop}[thm]{Proposition}

\theoremstyle{remark}
 
\newtheorem*{rem*}{Remark}
\theoremstyle{definition}

\newcommand{\todo}[1]{{\color{red}{#1}}}

\allowdisplaybreaks

\newtheorem{assumption}{Assumption}

\begin{document}
\title{Chance-Constrained Optimization for Non-Linear Network Flow Problems}
\author{Tillmann Weisser, Line A. Roald, and Sidhant Misra
\thanks{T. Weisser is with LAAS-CNRS, Toulouse. Email: tweisser@laas.fr}
\thanks{L. Roald and S. Misra are with Los Alamos National Laboratory. Email: \{roald,sidhant\}@lanl.gov .}}
\maketitle

\begin{abstract}
Many engineered systems, such as energy and transportation infrastructures, are networks governed by non-linear physical laws. A primary challenge for operators of these networks is to achieve optimal utilization while maintaining safety and feasibility, especially in the face of uncertainty regarding the system model. 
To address this problem, we formulate a Chance Constrained Optimal Physical Network Flow (CC-OPNF) problem that attempts to optimize the system while satisfying safety limits with a high probability. However, the non-linear equality constraints representing the network physics introduce modelling and optimization challenges which make the chance constraints numerically intractable in their original form. 
The main contribution of the paper is to present a method to obtain tractable polynomial approximations to the chance constraints using Semidefinite Programming (SDP). The method uses a combination of existing semi-algebraic techniques for projection and volume computation in combination with novel set manipulations to provide conservative inner approximations to the chance constraints. In addition, we develop a new two-step procedure to improve computational speed. 
While the method is applicable to general physical network flow problems with polynomial constraints, we use the AC optimal power flow problem for electric grids as an example to demonstrate the method numerically.
\end{abstract}


\section{Introduction}\label{sec:intro}

Networked systems are ubiquitous and include critical infrastructure networks such as the power grid, gas transmission pipelines, water networks and district heating systems. In such systems, optimization is often leveraged to maximize technical performance or economic efficiency, giving rise to what we will call Optimal Physical Network Flow (OPNF) problems.  
Optimization of system operation requires a mathematical model of the system. However, in practical systems, imperfect information and forecast errors introduce uncertainty in system operation and planning. 
If the uncertainty is not accounted for properly during the design and optimization process, the optimized system solution might be vulnerable to uncertainty, with potentially detrimental impacts on system risk. 

A prominent example is the Optimal Power Flow (OPF) problem in the electric power grid, which minimizes operational cost subject to technical constraints, and is used to clear electricity markets, perform security assessment and guide system expansion planning. The most significant source of uncertainty in the OPF problem is due to imperfect forecasts in renewable generation and loading conditions. System security must be maintained by ensuring that all variables are kept within acceptable values for a range of uncertainty realizations. Problems with similar structure also arise in other infrastructure networks such as natural gas and water networks.

A typical approach to account for uncertainty is to formulate the OPNF as a robust or stochastic program. However, for many of the above mentioned systems, the physics governing the network flows are given by a set of non-linear equations, such as branch flow equations and nodal conservation laws. This gives rise to non-linear equality constraints, which are inherently non-convex and thus challenging for both deterministic and stochastic optimization algorithms. In addition, the non-linearity significantly complicates the characterization of the uncertainty propagation throughout the system. 

Most existing robust and stochastic programming methods rely on assumptions of convexity. For practical problems such as the OPF problem, solution methods for robust or stochastic problem formulations typically use linear approximations \cite{roald2017, dallanese2017} or convex relaxations \cite{vrakopoulou2013AC, lorca2017robust, nasri2016} to circumvent the problem of non-convexity. This enables the application of well-known methods for robust \cite{robustoptimization2009} or chance-constrained \cite{campi2006, calafiore2006} programming, at the expense of a reduction in model fidelity and less comprehensive feasibility guarantees for the underlying problem. 

In this paper, we take a different approach. Instead of approximating or relaxing the non-linear network flow equations, we aim at treating the non-convex problem directly using techniques from polynomial optimization \cite{Henrion09, Magron15, Lasserre17}. The method is applicable for problems where the equality and inequality constraints can be represented as polynomials in both the decision variables and uncertain parameters.

We first formulate the uncertainty-aware problem as a Chance-Constrained OPNF (CC-OPNF) to guarantee that the constraints are satisfied with a high probability. 

Due to the chance-constraints, the CC-OPNF is intractable in the original form. The main contribution of the paper is to develop conservative, tractable approximations of the chance constraints in the form of polynomial constraints. We start from recent results for volume computations of semi-algebraic sets \cite{Henrion09} and projections of semi-algebraic sets \cite{Magron15}, which were shown to be useful for outer approximations of chance constraints \cite{Lasserre17}, and provide two crucial extensions: 
\begin{enumerate}
    \item While \cite{Lasserre17} allow for outer approximations of the chance constraints using a hierarchy of SDPs, in practice it is not straight-forward to obtain an \emph{inner approximation}, which is typically of interest in our setting.     Therefore, we use a series of set manipulations to extend the existing methods towards practical inner approximations. 
    \item To improve computational performance, we develop a \emph{two-step approximation procedure}, which allows for better approximations at lower computational overhead. 
\end{enumerate}

Replacing the chance constraints by their respective polynomial approximations yields an Approximate CC-OPNF (ACC-OPNF) problem that is still non-convex, but readily solvable by state-of-the-art non-linear programming solvers. Note that the polynomial chance-constraint approximations, which can be computationally heavy to compute, are determined a-priori in a \emph{pre-processing} step to the ACC-OPNF. 

Based on a small case study for the AC OPF problem, we demonstrate the practical performance of the method. In particular, we demonstrate the value of the extensions to inner approximations and the benefit of the two-step procedure. 

The remainder of the paper is organized as follows. After discussing the problem formulation in Section \ref{sec:probform}, we explain how to obtain polynomial approximations of the chance constraints in \ref{sec:PolApprox}. Section \ref{sec:improved_stokes} describes the rationale behind the two-step procedure, while Section \ref{sec:overall} summarizes the overall approach. In Section \ref{sec:application ccacopf}, we describe how our general framework can be mapped to the CC-AC-OPF, before providing numerical results in Section \ref{sec:numerics}. Finally, Section \ref{sec:concl} summarizes and concludes.


\section{Problem formulation}\label{sec:probform}
We now present the problem formulation in abstract form for a generic physical network flow problem, as the method can be applied to any problem that has the structure described below. 
For a concrete example, we refer the reader to section {\ref{sec:application ccacopf}}, where the method is applied to the AC~OPF problem. 

\subsection{Deterministic Optimal Physical Network Flow}
We define the problem variables as $\x = (x_1,\ldots,x_n)$ and $\y=(y_1,\ldots,y_m)$. For multivariate polynomials $f^0_1,\ldots,f^0_m$, $g^0_1,\ldots,g^0_k\in \rxy$ we consider the following Deterministic Optimal Physical Network Flow (D-OPNF) problem:
\begin{subequations}\label{eq:abstract OPF}
\begin{align}
    \min_{\x,\y}\;& c(\x,\y)\quad\st \nonumber \\
    & f^0_i(\x,\y) = 0 ,\;i=1,\ldots,m, \label{con:f)_deterministic} \\
    &  g^0_j(\x,\y)\geq0,\;j=1,\ldots,k. \label{con:g0_deterministic}
\end{align}
\end{subequations}
Here, the cost function is given by a polynomial $c\in \rxy$. The polynomial equality constraints $f^0_i(\x,\y)=0$ represent the network flow physics. The polynomial inequality constraints $g^0_j(\x,\y)\geq 0$ represent engineering limits.

To explicitly describe the degree of freedom in the system, we have separated the variables into $\x$ and $\y$. Since the equality constraints $f_i^0(\x,\y) = 0$ eliminate $m$ degrees of freedom, the variables $\y$ are an implicit function of the independent variable $\x$. Note that due to the non-linearity of the $f_i^0$, in general $\y$ might not be determined uniquely by a choice of $\x$. In this paper we make a practical assumption stated below.
\begin{assumption} \label{as:deterministic}
The engineering limits $g^0_j(\x,\y)\geq 0$ are such that the solution $\y$ to the system of equalities $f_i^0 = 0$, whenever it exists is unique. As a result, we therefore can write $\y$ as a function of $\x$, i.e. $\y_\x:=\y(\x)$.
\end{assumption}

The above assumption reflects a feature often encountered in engineered networks. Even though, mathematically the network physics described by the non-linear system $f^0_i = 0$ can have multiple solutions, there is only one solution that is physically meaningful within the region that the system is operated. As soon as the variables $\x$ are set, the state of the system is fully determined. Assumption~\ref{as:deterministic} allows us to formalize this notion.

\subsection{Chance-Constrained Optimal Physical Network Flow}
The aim of this paper is to account for uncertainty in the D-OPNF \eqref{eq:abstract OPF}, and to this end, we formulate the problem as a chance-constrained optimization problem. The chance constraints limit the probability of constraint violations, and can be enforced either as joint chance constraints (several equations hold jointly with a given probability) or separate chance constraints (each constraint is assigned its own probability). 
Due to the underlying physics of the problem, the network flow constraints $f^0_{i}$ must be satisfied jointly: If one of them is violated, the solution is not physically valid and the remaining constraints are meaningless. The probability of not jointly satisfying the network flow constraints can be understood as the probability that the uncertainty realization will lead to a situation where the flow problem is unstable and there exist no steady-state operating point (e.g., voltage instability in electric power grids).
The engineering limits $g_j^0$ can be satisfied either jointly or separately, depending on the preferred method for risk management. 
In this paper, we provide a method for enforcing the engineering limits as separate chance constraints. 

Let $(\Omega,\mu)$ be a probability space. The random variables $\w = (w_1,\ldots,w_\l)$ have zero mean $\boldsymbol{0}\in\R^l$. For every measurable event $A\subseteq\Omega$ denote the probability of $A$ by $\prob(A)=\int_A 1\dmu$. Finally, let $f_1,\ldots,f_m,g_1,\ldots,g_k\in \rxyw$ be multivariate polynomials. The notation is motivated by the idea that $f_i^0(\x,\y) = f_i(\x,\y,\boldsymbol{0})$ and $g_j^0(\x,\y) = g_j(\x,\y,\boldsymbol{0})$.
Define $f=\sum_{i=1}^m f_i^2$. Then enforcing the system of equations $f_i(\x,\y,\w) = 0, i=1,\ldots,m$ is equivalent to imposing $f=0$. We will use the later for better readability, although our implementation is based on the system of equalities rather than the single constraint $f=0$. We state the CC-OPNF problem:
\begin{subequations}\label{eq:abstract CCOPF interim}
\begin{align}
\min_{\x,\y_\x,\y(\w)}\;& c(\x,{\y_\x)}\quad\st \nonumber \\
&\hspace{-0.3cm}{f^0_i(\x,\y_\x) = 0 ,\;i=1,\ldots,m,} \label{con:f0_inter}\\
&\hspace{-0.3cm}{g^0_j(\x,\y_\x)\geq 0,\;j=1,\ldots,k,}\label{con:g0_inter}\\
&\hspace{-0.3cm}\prob\left(f(\x,\y(\w),\w) = 0\right)\geq 1-\varepsilon_1, \label{con:f_inter}\\
&\hspace{-0.3cm} \prob \left({g_j(\x,\y(\w),\w)\geq0}\right) \nonumber \\ 
    &\geq 1-\varepsilon_2,\;j=1,\ldots,k.\label{con:g_inter}
\end{align}
\end{subequations}
In addition to the chance-constraints to account for uncertainty, we also keep the constraints \eqref{con:f0_inter}, \eqref{con:g0_inter} from problem \eqref{eq:abstract OPF}. These constraints give a precise meaning to the cost function $c(\x,{\y_\x)}$, which is expressed as the operation cost for the expected realization $\w=\boldsymbol{0}$. 

We note that the problem as presented in \eqref{eq:abstract CCOPF interim} is a variational optimization problem. The $\x$ variables however do not depend on $\w$, which means that once the controllable variables are chosen, they cannot be modified in response to uncertainty. Although, the $\y$ variables are a function of $\w$, similar to the DNF \eqref{eq:abstract OPF}, the equality constraints eliminate the degrees of freedom for $\y$, and by a direct generalization of Assumption~\ref{as:deterministic}, one can think of $\y$ in \eqref{eq:abstract CCOPF} as a function of $(\x,\w)$, within the region defined by $g_j \geq 0$. As a result, the constraints in \eqref{con:g_inter} are simply constraints on the variable $\x$, a property that we exploit in our approach to convert the variational problem in \eqref{eq:abstract CCOPF} into a standard optimization problem in $\x$. However, attempting to eliminate the $\y(\w)$ variables creates another issue - unlike in \eqref{eq:abstract OPF}, where the inequalities \eqref{con:g0_deterministic} along with Assumption~\ref{as:deterministic} guarantee uniqueness and physical interpretability, eliminating $\y(\w)$ means that there is no way to enforce that $(\x,\y(\w),\w)$ satisfy all the inequalities in \eqref{con:g_inter}, thus forfeiting the aforementioned guarantees. 
To circumvent this issue, we introduce a set $Y$ for the $\y$ variables and make the following assumption:
\begin{assumption}
Restricting the range of $\y$ to a set $Y\subseteq\R^m$ the solution $\y(\w)$ to the system of equalities $f=0$ in \eqref{con:f_inter} is unique whenever it exists.
\end{assumption}
The set $Y$ can be interpreted as domain specific knowledge about the system introduced in order to reduce the feasible space to a region where our physical model is valid and exclude physically meaningless solutions to $f=0$. We propose the abstract formulation of the CC-OPNF problem below:
\begin{subequations}\label{eq:abstract CCOPF}
\begin{align}
\min_{\x,\y_\x}\;& c(\x,{\y_\x)}\quad\st \nonumber \\
&\hspace{-0.3cm}{f^0_i(\x,\y_\x) = 0 ,\;i=1,\ldots,m,} \label{con:f0}\\
&\hspace{-0.3cm}{g^0_j(\x,\y_\x)\geq 0,\;j=1,\ldots,k,}\label{con:g0}\\
&\hspace{-0.3cm}\prob\left({\exists \y\in Y,}\;f(\x,\y,\w) = 0\right)\geq 1-\varepsilon_1, \label{con:f}\\
&\hspace{-0.3cm} \prob \left({\exists\y\in Y,\;f(\x,\y,\w) = 0\land g_j(\x,\y,\w)\geq0}\right)\nonumber \\ 
    &\geq 1-\varepsilon_2,\;j=1,\ldots,k.\label{con:g}
\end{align}\end{subequations}
The sole reason for including the constraints $f=0$ in \eqref{con:g} is to implicitly specify $\y$ as a function of $\x$ and $\w$.

The main contribution of this paper is to provide \emph{tractable approximations} to the chance constraints \eqref{con:f}, \eqref{con:g}. The details of the approximation, which replaces the chance constraints in \eqref{con:f}, \eqref{con:g} by a set of polynomial constraints will be explained over the next sections.

\section{Polynomial approximations of chance constraints} \label{sec:PolApprox}
In this section, we first review results from the literature that use semi-definite programming (SDP) based methods for computing the volume of a basic semi-algebraic set, since they form the basis of the chance constraint approximations. Using these methods, we then develop inner and outer approximations of the chance constraint formulation in \eqref{eq:abstract CCOPF}.

\subsection{Preliminaries} \label{subsec:prelim}
Let $B = B_{\x} \times \Omega$, where $B_{\x} \subseteq \R^{n}$ is a hyper-interval or any other simple shape such that the moments with respect to the Lebesgue measure $\lambda_{\x}$ are known. We define a measure space $(B,\mu_{\x\w})$ by endowing $B$ with the product measure $\mu_{\x\w}$ given by $\mu_{\x\w} = \lambda_{\x} \otimes \mu$.  Let $K\subseteq B$ be a basic semi-algebraic set, where for each $(\x,\w) \in K$ we interpret $\x$ as the variables and $\w$ as the uncertainty. We call all points $(\x,\w) \in K$ as \emph{feasible} points. For a given $\x$, a chance constraint enforces that the probability that $(\x,\w)$ is feasible is larger than a given value, i.e.,
\begin{align}
    \prob((\x,\w) \in K) \geq 1-\epsilon, \label{eq:cc_general}
\end{align}
where the probability is computed using the measure $\mu$ on $\w$. This probability can be interpreted as the volume of the set $K_{\x} := \{\w : (\x,\w) \in K \}$ with respect to the measure $\mu$:
\begin{align}
    \rho(\x) := \prob((\x,\w) \in K) = \int_{\Omega} 1_{K_{\x}} d\mu,
    \label{eq:cc indicator}
\end{align}
where $1_{K_{\x}}$ denotes the indicator function of the set $K_{\x}$. 

\subsubsection{Approximating the volume of semi-algebraic sets}
In \cite{Henrion09} Henrion et al. propose a hierarchy of semi-definite programs approximating the set $K$ by the level set of some polynomial.
The starting point in \cite{Henrion09} is an infinite dimensional linear problem given as follows:
\begin{equation}\label{prob:full size dual}
\begin{split}
\min_{p\in\rxw} & \int_B p(\x,\w)\d\mu_{\x\w} \\
			\st & \quad p-1\geq 0 \mbox{ on } K,\\
				& \quad p \geq 0  \mbox{ on } B.
\end{split}
\end{equation}
Every feasible $p$ is an over-estimator of the indicator function of $K$ on $B$. By minimizing the integral over $p$, the optimal solution has to be close to the indicator function of $K$ in $L^1({\mu_{\x\w}})$-norm. The dual problem to \eqref{prob:full size dual} reads
\begin{equation} \label{prob:full size primal}
\begin{split}
\max_{\substack{\phi\in\meas(K)\\\psi\in\meas(B)}} & \int_K 1\d\phi \quad
\st\quad\forall (\alpha,\beta)\in \N_0^{n+\ell}\\
&\quad \int_K \x^\alpha\w^\beta \d\phi + \int_B \x^\alpha\w^\beta \d\psi = \int_B\x^\alpha\w^\beta \d\mu_{\x\w},
\end{split}
\end{equation}
where the optimization variables $\phi$ and $\psi$ are (positive) Borel measures supported on $K$ and $B$ respectively. As the moments, and in particular the mass of $\phi$ are bounded by the moments of ${\mu_{\x\w}}$, the optimal solution to \eqref{prob:full size primal} is the restriction of $\mu_{\x\w}$ to $K$. Consequently, the optimal value of \eqref{prob:full size primal} is the volume of $K$ with respect to $\mu_{\x\w}$. 

The infinite dimensional problems in \eqref{prob:full size primal} and \eqref{prob:full size dual} can be approximated by a hierarchy of semi-definite programs (SDPs) by using the method proposed by Lasserre \cite{Lasserre10}, which we briefly summarize.
A finite dimensional problem is obtained by (i) restricting the feasible set of \eqref{prob:full size dual} to polynomials of a fixed degree $2d$, and (ii) replacing the non negativity condition in the constraints by an algebraic certificate for non negativity (such as Putinar's theorem) on $K$ and $B$, respectively, which can be expressed by linear constraints on positive semi-definite matrices. The number $d$ is referred to as the relaxation degree or the relaxation order. The interested reader is referred to \cite{Lasserre10} for a full description of this relaxation procedure.

For any finite order of relaxation $d$ we obtain a polynomial $p_d\in\rxw$ of degree $2d$ that approximates the indicator function $1_K$ from above, and for any fixed $\x$ approximates the function $1_{K_{\x}}$ from above. Notice that since $p_d(\x,\w) \geq 1_K$ we have 
\begin{align} \label{eq:non-conservative}
   \rho(\x) =  \int  1_{K_{\x}} d \mu \leq \int p_d(\x,\w) \d \mu =: h^\ast(\x),
\end{align}
where the integration is only with respect to $\mu$, i.e, not with respect to $\mu_{\x\w}$. The chance constraint in \eqref{eq:cc_general} now can be replaced by the tractable polynomial inequality given by 
\begin{align} \label{eq:general_pover}
h^\ast(\x) \geq 1-\epsilon.
\end{align}
As $h^\ast$ is over approximating $\rho$, the constraint in \eqref{eq:general_pover} serves as an outer approximation of the chance constraints.

\subsubsection{Approximating the volume of the projection of semi-algebraic sets}
Comparing the generic representation of chance constraints in \eqref{eq:cc_general} to the one presented in \eqref{eq:abstract CCOPF}, we see that in many applications such as the OPNF, the presence of equality constraints introduce additional dependent variables $\y$ that are needed to describe the system. It is straightforward to extend the framework described above by appending the additional variables $\y$ to form the set $K$ in $(\x,\y,\w)$-space and apply the same procedure outlined above. However, since the variables $\y$ are fully specified by $(\x,\w)$ the volume of the set $K$ is zero leading to ill-conditioned problems while approximating the volume using \eqref{prob:full size primal}. This problem can be addressed by approximating the projection of $K$ onto the $(\x,\w)$ space where the volume is non-zero, instead of the original set $K$, using the method in Magron et al. \cite{Magron15}. To approximate the indicator function of the projection
\[
\pi_{\x\w}(K) :=\{(\x,\w) : \exists \y\in\R^m, (\x,\y,\w)\in K \}
\]
of $K$ onto the (\x,\w)-space, consider the variant of  \eqref{prob:full size dual}:
\begin{equation}\label{prob:projection dual}
\begin{split}
\min_{p\in\rxw} & \int_B p(\x,\w)\d\mu_{\x\y\w} \\
			\st & \quad p-1\geq 0 \mbox{ on } K,\\
				& \quad p \geq 0  \mbox{ on } B,
\end{split}
\end{equation}
where now $B:=B_\x\times B_\y\times\Omega$ for some set $B_\x$ and $B_\y$ for which it is easy to compute the moments of the Lebesgue measure, and $\mu_{\x\y\w} = \lambda_{\x} \otimes \lambda_{\y} \otimes \mu$.
Note that the optimizing variable $p$ is restricted to be invariant in $\y$-direction. The constraints guarantee that $p$ is an over-estimator of the indicator function of $\pi_{\x\w}(K)$ on $\pi_{\x\w}(B)=B_\x\times\Omega$. Similar to the results in {\cite{Henrion09}}, Magron et al. prove convergence results for $p$ to the indicator function of $\pi_{\x\w}(K)$  and the optimal value to the volume of the projection with respect to the marginal of $\mu_{\x\y\w}$ for the corresponding semi-definite hierarchies. 


\subsection{Approximations of the CC-OPNF}
In this subsection, we describe how to use the methods outlined in Section~\ref{subsec:prelim} to provide outer and inner approximations of the CC-OPNF \eqref{eq:abstract CCOPF}.
We specify the feasible set of the chance constraints that we want to approximate by 
\begin{subequations}\label{eq:formulation CCOPF}
\begin{align}
\kxp&:=\{\x: \prob(\exists \y\in Y, f(\x,\y,\w) = 0)\geq 1-\varepsilon_1, \label{lx:1}\\ 
&\prob (\exists \y\in Y, f(\x,\y,\w) = 0\land g_j(\x,\y,\w)\geq0)\nonumber\\
    &\quad \geq 1-\varepsilon_2,\quad j=1,\ldots,k\}, \label{lx:2}
\end{align}
\end{subequations}
where we assume that $\kxp\subseteq B_\x$. As mentioned in Section~\ref{sec:probform}, our goal is to approximate the set $\kxp$ by replacing the intractable chance constraints by polynomial constraints. We define the sets for which the constraints remain satisfied as
\begin{align}\label{eq:kj outer}
K_0&:=\{(\x,\y,\w)\in B: f(\x,\y,\w) = 0 \},\\
K_j&:=\{(\x,\y,\w)\in B: f(\x,\y,\w) = 0\land g_j(\x,\y,\w)\geq 0 \}, \nonumber \\
    & \qquad j=1,\ldots,k. \nonumber 
\end{align}

\subsubsection{Outer approximation of the feasible set}\label{sec:outer}
An outer approximation of the set $\kxp$ can be obtained by applying the method outlined in Section~\ref{subsec:prelim} to each of the sets $K_j$ for $j=0,\ldots,k$. For each $K_j$, we get a polynomial $h^\ast_j\in\rx$ which approximates the function $\x\mapsto\prob(\pi_{\x\w}(K_j))$ from above, leading to an overestimation of the satisfaction probability and an outer approximation of the chance constraints. Consequently the set 
\[
\{x\in B_\x: h^\ast_0(\x)\geq1-\varepsilon_1,h^\ast_j(\x)\geq1-\varepsilon_2,j=1,\ldots,k\}
\]
is an outer approximation of \kxp, and the corresponding ACC-OPNF provides a lower bound to the optimal cost of the OPNF.

\subsubsection{Inner approximation of the feasible set}\label{sec:inner}
In applications where system security is of primary concern, obtaining feasible solutions to \eqref{eq:abstract CCOPF} are more important than obtaining lower bounds to the cost, motivating an investigation of inner approximations to the chance constraints. However, as opposed to the outer approximation, obtaining an inner approximation of \kxp\ is more involved. 

In the following, we propose a modification of \kxp\, which we can use to approximate (almost) from the interior.
For $\varepsilon_1<\varepsilon_2$ define the set 
\begin{subequations}
\begin{align}
\kx& := \{\x\in B_\x : \nonumber\\ 
&\prob(\exists \y\in Y, f(\x,\y,\w)= 0) \geq 1-\varepsilon_1, \label{kx:1}\\ 
&\prob(\exists \y\in Y, (f(\x,\y,\w)=0 \land g_j(\x,\y,\w) \leq 0)) \nonumber\\
&\quad \leq \varepsilon_2-\varepsilon_1,\quad j=1,\ldots,k\}. \label{kx:2}
\end{align}
\end{subequations}

The essential difference between \kxp\ and \kx\ is that the probabilities \eqref{kx:2} in \kx\ are bounded from above whereas the probabilites \eqref{lx:2} in \kxp\ are bounded from below. Since the methods discussed in \ref{subsec:prelim} lead to over-estimators of the probability, the reversal of the inequality in the formulation in \kx\ now enables us to approximate the sets described by the chance constraints in \eqref{kx:2} from the interior.
The following proposition relates the approximating set \kx\ to \kxp.

\begin{prop}
\kx\ is an inner approximation of \kxp.
\end{prop}

The proof is simple and is given in the appendix. Instead of directly dealing with $\kxp$, we attempt to approximate the set $\kx$ from the interior.
Using the same procedure we now compute polynomials $h^\ast_0,\ldots,h^\ast_m$ approximating the functions $\x\mapsto\prob(\pi_{\x\w}(K_j))$ where $K_j$ now is defined by

\begin{align}\label{eq:kj inner}
\kxyw_0:=& \{(\x,\y,\w)\in B:f(\x,\y,\w)=0 \},\\
\kxyw_j:=& \{(\x,\y,\w)\in B:f(\x,\y,\w)=0 \land g_j(\x,\y,\w)\leq 0 \}.\nonumber
\end{align}

Note that though we are aiming for an inner approximation of \kx, the polynomials $h^\ast_j$ are over-approximators of the probability. The set \kx\ is then approximated by the set
\begin{subequations}
\begin{align}
\tilde\kx := \{ x \in B_\x:  &\; h^\ast_0(\x)\geq 1-\varepsilon_1 \label{eq:equality_approx}\\
 &\; h^\ast_j(\x)\leq \varepsilon_2-\varepsilon_1,j=1,\ldots,k \label{eq:inequality_approx}\}.
\end{align}
\end{subequations}
Since the polynomials $h^\ast_j(\x)$ over-approximate the probabilities in \eqref{kx:2}, the set defined by the inequalities in \eqref{eq:inequality_approx} are inner approximations of the corresponding sets defined by \eqref{kx:2}. Unfortunately the same relation is not true for the sets defined by \eqref{eq:equality_approx} and \eqref{kx:1} that correspond to the probability of joint violation of the equatility constraints $f_i(\x,\y,\w) = 0$. Therefore, $\tilde\kx$ is an \emph{approximate} inner approximation to \kxp.

\section{Improved approximations through Stokes constraints}
\label{sec:improved_stokes}
The SDP hierarchy to approximate the chance constraints presented in Section~\ref{subsec:prelim} is guaranteed to converge to the optimum as $d$ grows to infinity, but much less is known about the associated rate of convergence. When the number of variables in the polynomial optimization problem is large, the computation times can become prohibitively expensive, since current SDP-solvers are not able to solve problems with variables of size $>1000$ on a standard computer. {Coupled with the fact that the size of SDP-variables at relaxation level $d$ is $\binom{N+d}{d}$,where $N$ is the number of variables of the polynomial optimization problem, it becomes crucial to achieve high approximation accuracy at lower values of $d$.}
However, convergence of the indicator function approximation tends to be slow due to the so-called Gibbs' phenomenon:  
If a function has a discontinuity, every overestimating polynomial approximation $p$ overshoots the upper value at the jump \cite{Lasserre17}. In the following subsections, we first review existing results regarding the use of valid constraints generated via the Stokes Theorem to speed up the convergence rate, and then describe our approach to generalize this procedure to computing the volume/probability of projections of semi-algebraic sets. 

\subsection{Concept of Stokes constraints}\label{subsec:stokes for cc}
In \cite{Lasserre17} Lasserre proposes to improve convergence of the hierarchy by adding additional constraints to the problem \eqref{prob:full size primal}. When a polynomial $t$ is known to vanish on the boundary of $K$, the optimal measure $\phi^\ast$ satisfies the equality $\int_K\theta(\x,\y)\d\phi^\ast = 0$ for some family of functions $\theta$ depending on $\mu$ and $t$. The equality is a consequence of the Stokes theorem,  which is why the constraints are referred to as Stokes constraints. {We describe the procedure to generate these constraints below, in the case where $\mu$ is the uniform measure. For more general measures we refer to \cite{Lasserre17}.}

Let $t\in \R[\x,\w]$ be a polynomial that vanishes on the boundary of $K$. Given any $(\alpha,\beta) \in \N_0^{n+\ell}$ and $z \in\{x_1,\ldots,x_n,w_1,\ldots,w_m\}$, define the polynomial $\theta_{\alpha\beta}^z$ as
\begin{align} \label{eq:theta_def}
    \theta_{\alpha\beta}^z := \tfrac{\partial}{\partial z}\left(\x^\alpha\w^\beta t(\x,\w) \right).
\end{align}
Then by the Stokes formula, for all $(\alpha,\beta) \in \N_0^{n+\ell}$ and $z \in\{x_1,\ldots,x_n,w_1,\ldots,w_m\}$ we have 
\begin{align} \label{prob:stokes primal}
    \int_K \theta_{\alpha,\beta}^z \d \phi^\ast= \int_{ K}\tfrac{\partial}{\partial z}\left(\x^\alpha\w^\beta t(\x,\w)\right) \d \phi^\ast = 0. 
\end{align}
Since the optimal measure satisfies all the equality constraints given in \eqref{prob:stokes primal}, we can add these equations as constraints to \eqref{prob:full size dual} without affecting the optimal solution. Adding these redundant constraints has been shown in some cases to greatly improve the rate of convergence of the SDP hierarchy, i.e., enabling higher accuracy at a lower relaxation level $d$. 
While the faster convergence is beneficial, the dual of \eqref{prob:full size dual} with addition of constraints in \eqref{prob:stokes primal} now reads
\begin{equation}\label{prob:stokes dual}
\begin{split}
\min_{\substack{p\in\rxw,\\q_{\theta}\in \R}} & \int_B p(\x,\w)\d\mu_{\x\w} \\
			\st & \quad p-1\geq \sum_{\theta\in \Theta}q_\theta\theta \mbox{ on } K,\\
				& \quad p \geq 0  \mbox{ on } B,
\end{split}
\end{equation}
where $\Theta$ is the set of all $\theta^z_{\alpha\beta}$ defined in \eqref{eq:theta_def}. Comparing \eqref{prob:stokes dual} to problem \eqref{prob:full size dual} we observe that the polynomial $p$ in \eqref{prob:stokes dual} is no longer an over-estimator of the indicator function $1_K$ on $B$. 

\begin{figure}
\begin{center}
\includegraphics[width = 0.48\textwidth]{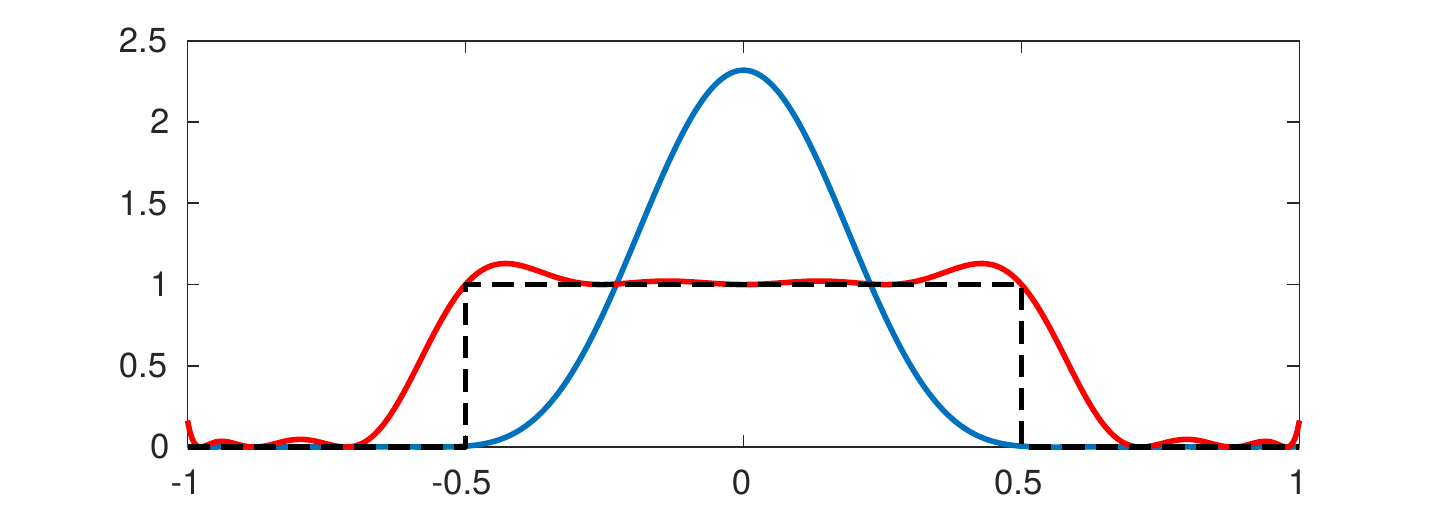}
\caption{Effect of Stokes constraints on the dual variables. Polynomial $p$ approximating the indicator function on $[-\tfrac{1}{2},\tfrac{1}{2}]$ without (red) and with (blue) Stokes constraints. }
\label{fig:stokes effect}
\vspace{-12pt}
\end{center}
\end{figure}

This effect is illustrated in Figure~\ref{fig:stokes effect}, where typical shapes of $p$ for problems \eqref{prob:full size dual} (red) and \eqref{prob:stokes dual} (blue) are shown. We observe that the red curve over-approximates the indicator function of $K$ (dashed black). Also we can see the mismatches at the discontinuities of the indicator function due to the Gibbs' phenomenon. We note that the $1$-super-level-set of the red curve is a good approximation of the set $K$. In contrast to that, when applying Stokes constraints, we observe that the $1$-super-level-set of the blue curve does not provide any information about the set $K$. The integral value of the blue curve however is closer to the volume of $K$ than the integral value of the red curve. Moreover the integral preserves the over-approximation property. 
Indeed for any polynomial $p$ feasible for \eqref{prob:stokes dual} we have
\begin{align} \label{eq:stokes_volume_upper}
    \int_K p \d\mu_{\x\w} \geq \int_K 1 \d\mu_{\x\w} +\sum_{\theta\in \Theta}q_\theta\hspace{-6pt}\underbrace{\int_K\hspace{-6pt} \theta\d\mu_{\x\w}}_{=0\;\text{by Stokes}}\hspace{-6pt} {=} \vol(K).
\end{align}

\subsection{Partial Stokes constraints for chance constraints} \label{sec:partial_stokes}
The polynomial $p(\x,\w)$ obtained above cannot be used to approximate the function $\rho(\x)$ in \eqref{eq:cc indicator}.
If however we only use $\w$ as variables $z$ in \eqref{eq:theta_def} then for all $\x$ we have 
\begin{align} \label{eq:partial_stokes}
    \int_{B_{\x}} p d\mu \geq \int_{K_{\x}} 1 \d \mu + \sum_{\theta \in \Theta} \int_{K_{\x}} \theta \d \mu = \rho(\x).
\end{align}
Applying Stokes constraints only in the $\w$ direction hence allows us to both obtain the improved convergence rates while still obtaining an over-estimator of the probability of $K_{\x}$.

\subsection{Partial Stokes constraints for projection of sets}\label{sec:two steps}
The method in Section~\ref{subsec:stokes for cc} cannot directly be applied to the setting where the feasible set is described by the projection of a semi-algebraic set. This is because in order to be able to add Stokes constraints to the problem in \eqref{prob:projection dual}, we must first find a polynomial $t\in\rxw$ that vanishes on the boundary of the projection $\pi_{\x\w}(K)$ of $K$.
Note that in Section~\ref{subsec:stokes for cc} where there is no projection involved, the polynomial $t\in\rxw$ that vanishes on the boundary of $K$ can be readily obtained as the product of the polynomials that define the semi-algebraic set $K$.
For the projection $\pi_{\x\w}(K)$ of a semi-algebraic set $K$ in $(\x,\y,\w)$-space, this trick is not applicable.

Our solution to this issue is a two-step-procedure: In the first step we approximate the projection $\pi_{\x\w}(K)$ by the super-level-set of a polynomial $p^{(1)}\in\rxw$. In the second step, we use this super-level-set $S$ to compute a second polynomial $p^{(2)}\in\rxw$ approximating the volume of $S$. This is explained in more detail below.

\subsubsection{Step 1: Approximating the projection of $K$}\label{sec:step1}
\begin{figure}	
\begin{center}
		\def\svgwidth{0.4\textwidth}
		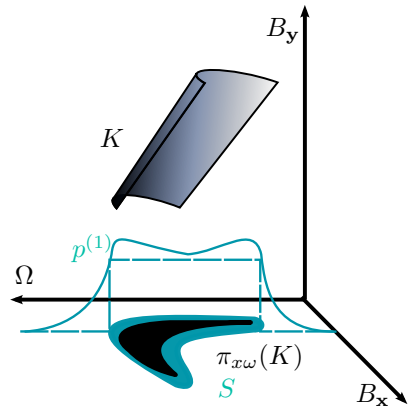
		\caption{Step 1: Projection step. The projection of $K$ is approximated as $S$, which is defined by the 1-super-level set of $p^{(1)}$.}
		\label{fig:sketch1}
\end{center}
\end{figure}

We first apply the method in Section~\ref{subsec:prelim} and solve the problem in \eqref{prob:projection dual} to obtain a polynomial $p^{(1)}$ that is an over-estimator of the indicator function of $\pi_{\x\w}(K)$, i.e. $p^{(1)}\geq1 \text{ on } \pi_{\x\w}(K)$.
In particular the super-level-set given by
\begin{align}
    S:=\{(\x,\w)\in B_\x \times \Omega : p^{(1)}(\x,\w)-1\geq0 \}   \label{eq:S_def}
\end{align}
is an outer approximation of $\pi_{\x\w}(K)$. Fig.~\ref{fig:sketch1} illustrates this step. Numerical experiments have shown that the $1$-super-level-set of the optimizing polynomial is quite accurate already for low relaxation degrees.

\subsubsection{Step 2: Probability approximation}\label{sec:step2}
\begin{figure}	
\begin{center}
		\def\svgwidth{0.4\textwidth}
		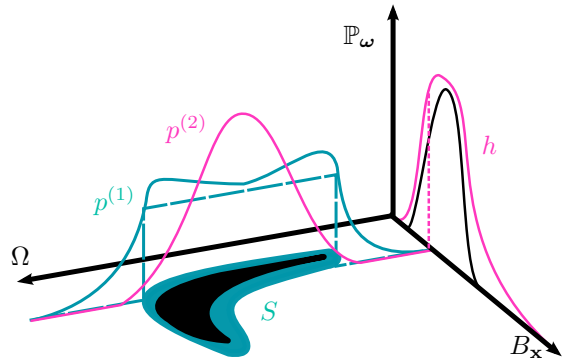
		\caption{Step 2: Probability approximation. The probability is approximated by integrating $p^{(2)}$ in $\Omega$ direction for every $\x$.}
		\label{fig:sketch2}
\end{center}
\end{figure}

After the first step we replace the actual projection $\pi_{\x\w}(K)$ by its approximation $S$ defined in \eqref{eq:S_def}. In doing so we lose information about $\pi_{\x\w}(K)$ but we gain two important advantages. First, moving from $K$ to $S$ we get a significant reduction in the number of variables, as we eliminate the whole $\y$-space. This allows us to afford computational capacity for higher levels in the SDP relaxation hierarchy and get better volume approximations. Second, we now have a polynomial, specifically $p^{(1)}-1$, that vanishes on the boundary of $S$. This crucial difference enables us to use Stokes constraints to improve the volume approximation. Applying the method in Section~\ref{sec:partial_stokes}, we obtain a polynomial $p^{(2)}\in\rxw$ that still preserves the desired over-approximation property:
\begin{align}
    h^\ast(\x):=\int_{\Omega} p^{(2)}(\x,\w)\d\mu \stackrel{(a)}{\geq} \prob(S) \stackrel{(b)}{\geq} \prob(\pi_{\x\w}(K)), \nonumber
\end{align}
\todo{}
where $(a)$ follows from \eqref{eq:partial_stokes} and $(b)$ follows because ${\pi_{\x\w}(K)} \subseteq S$. 
This step is summarized in Fig.~\ref{fig:sketch2}.


\section{The overall approach}\label{sec:overall}
To summarize the overall approach, we first recall the problem formulation \eqref{eq:abstract CCOPF}. Our aim is to eliminate the chance constraints {\eqref{con:f} and \eqref{con:g}} and replace them by tractable polynomial constraints. The challenge is to (i) ensure existence of solution to the equality constraints, (ii) compute inner approximations to the chance constraints, and (iii) enable use of Stokes constraints to speed up convergence.

We address the challenges in the following steps:
\begin{enumerate}
    \item We reformulate the feasible set $\kxp$ of the chance constraints by the set \kx\ that allows us to obtain inner approximations.
    \item We eliminate the dependent $\y$ variables by approximating the projection of each $K_j$ defining \kx\ as the super-level set $S$ of a polynomial $p_j^{(1)}$.
    \item We use the reduced set $S$ to compute the inner approximations to the chance constraints by polynomials $h^\ast_0(\x),\ldots,h^\ast_k(\x)$. To speed up convergence, we add Stokes constraints which is made possible by the availability of the polynomial $p^{(1)}$. 
\end{enumerate}

Now, the chance constraints in original problem \eqref{con:f}, \eqref{con:g} can be replaced by their approximation to obtain the ACC-OPNF formulation:
\begin{subequations}\label{eq:acc-opnf}
\begin{align}
\min_{\x,\y_\x}\;& c(\x,{\y_\x)}\quad\st \nonumber \\
&{f^0_i(\x,\y_\x) = 0 ,\;i=1,\ldots,m,}\\
&{g^0_j(\x,\y_\x)\geq 0,\;j=1,\ldots,k,}\\
& h_0(\x) \geq 1-\varepsilon_1,\label{con:h0}\\
& h_j(\x) \leq \varepsilon_2-\varepsilon_1,\; j=1,\ldots,k.\label{con:hj}
\end{align}\end{subequations}

Although obtaining the polynomials $h^\ast_0(\x),\ldots,h^\ast_k(\x)$ might be computationally heavy, this procedure is independent of the actual solution process for the resulting ACC-OPNF and can be considered as a pre-processing step to be executed offline.  
The resulting approximate CC-OPNF, despite remaining non-convex, can be solved to local optimality easily using a local non-linear solver. Furthermore, methods for global optimization of polynomial problems can be applied \cite{Lasserre2001}.

\section{Application to Chance-Constrained AC Optimal Power Flow}\label{sec:application ccacopf}
In this section, we present the mapping of a chance-constrained AC optimal power flow (CC-AC-OPF) problem onto the general CC-OPNF problem \eqref{eq:abstract CCOPF}. 
Motivated by the recent increase in generation uncertainty from renewable energy sources, our CC-AC-OPF formulation attempts to minimize generation cost, subject to engineering constraints while accounting for the uncertainty in  renewable power generation. 

\subsection{Deterministic Optimal Power Flow} 
We first formulate the deterministic OPF problem where we assume perfect knowledge of the system. This problem corresponds to the deterministic OPNF \eqref{eq:abstract OPF}. 

\subsubsection{Notation} We consider an electric network where $\mathcal{N}$ and $\mathcal{E}$ denote the sets of nodes and edges. 
Without loss of generality, we assume that there is one generator, one demand and one uncertainty source per bus.
Complex power is given by $s=p+j\cdot q$, where $p$ and $q$ are the active and reactive power. Subscripts $_R,~_G$ and $_D$ are for renewable energy sources, conventional generators and loads, respectively. The complex bus voltages are denoted by $v=v_{real}+v_{imag}$, and the corresponding voltage magnitudes by 
$|v|=(v_{real}^2+v_{imag}^2)^{1/2}$.

\subsubsection{Problem formulation} Given the above considerations, the OPF problem is given by
\begin{subequations}
\begin{align}
     \min_{\substack{p_{G0}, \\ q_{G0},v_0}} ~&\sum_{i\in\mathcal{G}} c_{2,i} p_{G0,i}^2 + c_{1,i} p_{G0,i} + c_0 \label{eq:opfobjective}\\
     \text{s.t.} ~~~
    &s_{G0,i}+ s_{R,i} - s_{D,i} = \sum_{(i,j)\in\mathcal{E}} s_{0,ij}, && \forall i\in\mathcal{N}, \label{eq:acnodal}\\
    &s_{0,ij} = \Y_{ij}^* v_{0,i} v_{0,i}^* - \Y_{ij}^* v_{0,i} v_{0,j}^*, && \forall (i,j) \!\in \!\mathcal{E}, \label{eq:acflow}\\
    &p_{G,i}^{min} \leq p_{G0,i} \leq p_{G,i}^{max}, \quad  &&\forall i\in\mathcal{N}, \label{eq:pG}\\
    &q_{G,i}^{min} \leq q_{G0,i} \leq q_{G,i}^{max}, \quad &&\forall i\in\mathcal{N}, \label{eq:qG}\\
    &|v|^{min} \leq |v_{0,j}| \leq |v|^{max}, \quad &&\forall j\in\mathcal{N}, \label{eq:v}\\
    &|s_{0,ij}| \leq |s_{ij}|^{max}, \quad &&\forall (i,j)\!\in\!\mathcal{E}. \label{eq:s0}
\end{align}
\label{eq:detOPF}
\end{subequations}
The objective \eqref{eq:opfobjective} of the problem is to chose the generation dispatch point, given by the active and reactive power generation $p_{G0},~q_{G0}$ and the complex voltages $v_0$, such that the the cost of active power generation given by the quadratic function in \eqref{eq:opfobjective} is minimized.  
The AC power flow equations \eqref{eq:acnodal}, \eqref{eq:acflow} are a set of equality constraints describing the physical laws, with the nodal power balance given by \eqref{eq:acnodal}, and transmission line flows given by the Ohm's law \eqref{eq:acflow}, where $\Y$ is the so-called \emph{admittance matrix}. Note that we use the rectangular form of the power flow equations to obtain polynomial constraints. 
Further, we enforce a set of engineering limits \eqref{eq:pG}-\eqref{eq:s0}. The constraints \eqref{eq:pG}, \eqref{eq:qG} represent bounds on generation capacity, \eqref{eq:v} limits the voltage magnitudes to safe ranges and \eqref{eq:s0} enforces limits on the apparent power flow. Among thesse constraints, \eqref{eq:acnodal} and \eqref{eq:acflow} correspond to the equality constraints $f_i^0=0$ in the deterministic OPNF \eqref{eq:abstract OPF}, and the remaining constraints correspond to the inequality constraints $g^0_j\geq 0$.

\subsection{Chance-Constrained Optimal Power Flow}
We now extend the deterministic problem to the setting with uncertainty in the power injections.
\subsubsection{Modelling uncertain injections} 
We model the uncertain active power injections from renewable generators as the sum of the expected value $p_R$ and a fluctuation $\w$. The expected reactive power injection is denoted by $q_R$.  The reactive power injections are assumed to adjust in a way that the power factor, given by $\gamma=q_R/p_R$, remains constant:
\begin{equation}
    s_R(\w) = (p_R + \w) + j \cdot (q_R + \gamma \w) \label{eq:power_factor} 
\end{equation}
We assume that the probability distribution of $\w$ is known. The active and reactive power consumption of the loads, denoted by $p_L,~q_L$, are assumed to be constant, but could also be modelled similar to \eqref{eq:power_factor}.

\subsubsection{Power flow equations under uncertainty}
For non-zero uncertainty realization $\w$, the power flow equations \eqref{eq:acnodal} are adapted to account for $\w$, i.e.
\begin{subequations}
\label{eq:powerflow}
\begin{align}
    &\!\!\!s_{G,i}(\w) + s_{R,i} + \w - s_{D,i} \!= \!\!\sum_{(i,j)\in\mathcal{E}} s_{ij}(\w), \!\!&& \!\forall i\in\mathcal{N}, \\
    &\!\!\!s_{ij}(\w) \!=\! \Y_{ij}^* v_i(\w) v_i^*(\w) \!-\! \Y_{ij}^* v_i(\w) v_j^*(\w), \!\!&& \!\forall (i,j)\! \in\! \mathcal{E}. 
\end{align}
\end{subequations}

\subsubsection{Response to uncertainty} 
When the power injections fluctuate, the controllable generators must adjust their generation output $s_{G,i}(\w)$ to ensure that the power balance constraints \eqref{eq:acnodal} are satisfied. 
We adopt balancing practices typical in power systems operation, which require the definition of so-called $pv$, $pq$ and $v\theta$ (reference) buses.

On each node of the network there are four state variables, namely the active power injection $p$, the reactive power injection $q$, and two voltage variables corresponding to the voltage magnitude and angle $|v|,~\theta$ (polar coordinates) or the real and imaginary voltage $v_{real},~v_{imag}$ (rectangular coordinates). The buses are classified according to the quantities that are controllable or specified: (i) $pq$ buses (such as loads) with specified real and reactive power, (ii) $pv$ buses (such as generators) with controllable active power and voltage magnitude, and (iii) $v\theta$ or reference bus with the voltage angle set to zero. The sets of nodes that correspond to the three categories are denoted by subscripts $\mathcal{N}_{pq},~\mathcal{N}_{pv}$ and $\mathcal{N}_{v\theta}$.

Given the above definitions, we assume that the active power injections from generators at $pq,~pv$ buses remain constant throughout the fluctuations, and all fluctuations $\omega$ are balanced by the generator connected at the slack bus. 
Similarly, reactive power is balanced by adjusting the reactive power output of $pv$ and $v\theta$ buses to maintain constant voltage magnitudes, while the reactive power injections at $pq$ buses are kept constant.

\subsection{Definition of $\x$ and $\y$ variables}
We choose the rectangular coordinate representation in order to be able to employ the semi-algebraic methods described in this paper. This gives us $4$ variables per bus $p,q,v_{imag},v_{real}$. However, as described above, the standard model for $pv$ and $v\theta$ buses are based on polar coordinates, where we keep the voltage magnitude constant. We handle these requirements in rectangular coordinates by adding the constraints $v_{imag}=0$ and $v_{real,i}(\w)=v_{real,i}$ for $i\in\mathcal{N}_{v\theta}$, and the constraint $v_{real,i}(\w)^2 + v_{imag,i}(\w)^2 = |v|_i^2$ for $i\in\mathcal{N}_{pv}$.

This results in two independent variables per bus, which we choose to also correspond to the quantities that can be controlled by the system operator. In particular, we define the independent $\x$ variables as 
\begin{align*}
    & p_{G0,i}, q_{G0,i}, ~ &&\forall i\in\mathcal{N}_{pq}, \\
    & p_{G0,i}(\w),|v|_{0,i} , ~ &&\forall i\in\mathcal{N}_{pv},\\
    & v_{\text{real}0,i}, v_{\text{imag}0,i},  ~ &&\forall i\in\mathcal{N}_{v\theta}.
\end{align*}
The variables that change as a function of $\w$ are the $\y$ variables in the CC-OPNF formulation \eqref{eq:abstract CCOPF}: 
\begin{align*}
    & v_{\text{real},i}(\w), v_{\text{imag},i}(\w) , ~ &&\forall i\in\mathcal{N}_{pq}, \\
    & q_{G,i}(\w),v_{\text{real},i}(\w),v_{\text{imag},i}(\w) , ~ &&\forall i\in\mathcal{N}_{pv},\\
    & p_{G,i}(\w),q_{G,i}(\w), ~ &&\forall i\in\mathcal{N}_{v\theta},\\
    & s_{ij}(\w), ~ && \forall ij\in\mathcal{E}.
\end{align*}
Note that in the process of solving \eqref{eq:abstract CCOPF}, we are not explicitly assigning a value to these dependent quantities $\y(\w)$.
However, the variables $\y_\x$, which correspond to the $\y$ variables at the expected operating point ($\w=0$), are explicitly defined.

\subsubsection{Definition of constraints $f=0$ and $g\leq0$}
As is evident from \eqref{eq:powerflow}, both the generation outputs $p_{G,i}(\w)$ and $q_{G,i}(\w)$, the power flows $s_{ij}(\w)$ and the voltage variables $v_i(\w)$ will change depending on the realization of $\w$. The constraints which incorporate those quantities are therefore enforced as chance constraints. 

The stochastic power flow equations \eqref{eq:powerflow} correspond to the equality constraints $f(\x,\y,\w) = 0$.
When there is no solution to this set of equations, the system is unstable and might collapse at any point leading to complete blackout of the electric grid. We hence want the probability of violating any of the equality constraints to be very low, and enforce those constraints jointly as in \eqref{con:f} with a small acceptable violation probability $\varepsilon_1$. 

The inequality constraints $g_j(\x,\y,\w)\leq 0$ correspond to the engineering limits
\begin{subequations} 
\label{eq:engineeringCC}
\begin{align}
    &p_{G,i}^{min} \leq p_{G,i}(\w)  \leq p_{G,i}^{max}, \quad &&\forall i\in\mathcal{N}_{v\theta} \label{eq:pGref}\\
    &q_{G,i}^{min} \leq q_{G,i}(\w)  \leq q_{G,i}^{max}, \quad &&\forall i\in\mathcal{N}_{pv}, \mathcal{N}_{v\theta} \label{eq:qGpv}\\
    &|v_i|^{min} \leq |v_{i}|(\w)  \leq |v_i|^{max}, \quad &&\forall i\in\mathcal{N}_{pq} \label{eq:vpq}\\
    &v_{real,i}(\w)^2 + v_{imag,i}(\w)^2 = |v|_i^2,  \quad &&\forall i\in\mathcal{N}_{pv} \label{eq:vpv}\\
    &|s_{ij}|(\w)  \leq |s_{ij}|^{max}, \quad &&\forall (i,j)\in\mathcal{E}. \label{eq:s}
\end{align}
\end{subequations}
In contrast to a violation of the power flow equations \eqref{eq:powerflow}, a violation of one of the engineering constraints \eqref{eq:engineeringCC} would typically have a more local impact (e.g. overloading of a component), and can often be tolerated for a certain amount of time (e.g. violations of thermal capacity limits of transmission lines). We hence enforce \eqref{eq:engineeringCC} as separate chance constraints, and allow for a larger violation probability  $\epsilon_2 > \epsilon_1$.

\subsubsection{Choosing $Y$}
The last parameter we must determine before the mapping from the CC-AC-OPF to the generic CC-OPNF problem \eqref{eq:abstract CCOPF} is complete, is the set $Y$ from Assumption 2. We would like to choose $Y$ such that solutions to \eqref{eq:powerflow} are unique and have a well-defined physical meaning, which for the OPF problem implies ensuring that low voltage solutions to the power flow equations are excluded. Therefore we define the sets $Y$ by the inequalities
\begin{equation}\label{eq:anti low voltage}
    |v|^{min-} \leq |v_{i}|(\w),\quad\forall i\in\mathcal{N}_{pq}.
\end{equation}
Here, $|v|^{min-}$ is lower than the standard voltage bound $|v|^{min}$, but sufficiently large to exclude low voltage solutions.

\section{Case study}\label{sec:numerics}
We first describe the implementation and test system, before presenting the numerical results for the chance constraint approximation and the resulting approximate CC-OPNF.

\subsection{Implementation}
In this section, we describe our implementation to obtain the ACC-OPNF in Section~\ref{sec:overall} and evaluate its performance. To obtain the polynomials $h^\ast_0,\ldots,h^\ast_k$ in \eqref{eq:acc-opnf} we solve SDP relaxations to the infinite dimensional linear problems described in \ref{sec:step1} and \ref{sec:step2}. We use the GloptiPoly3 Matlab toolbox \cite{gloptipoly} to model the relaxations and Mosek \cite{mosek} to solve the SDPs.
The resulting ACC-OPNF is implemented in Julia \cite{julia} with JuMP \cite{JuMP} and PowerModels.jl \cite{PowerModels} and then solved using the local non-linear solver Ipopt \cite{ipopt}. We also perform Monte-Carlo simulations for benchmarking which requires solving the standard power flow and the AC-OPF which are implemented using Matpower \cite{zimmermann2011} and PowerModels.jl respectively.

\subsection{Test system}
\begin{figure}
\centering
\includegraphics[width = 0.9\columnwidth]{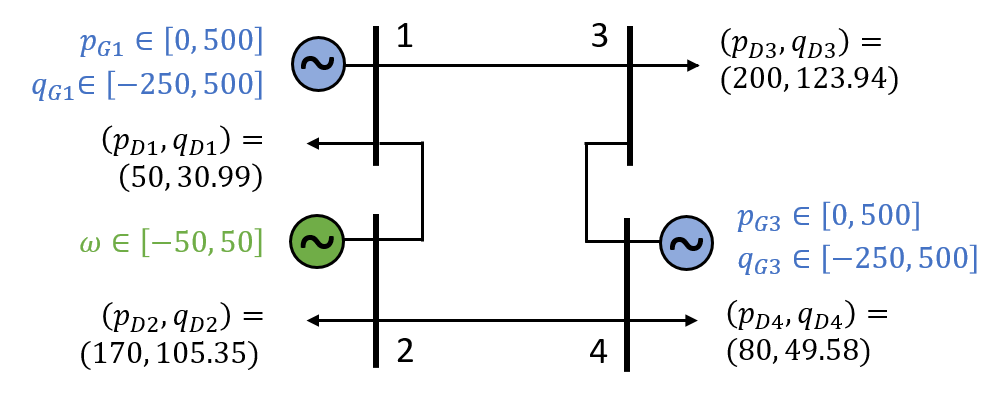}
\caption{Overview of the $4$-bus system. Generators marked in blue, uncertainty source in green and loads in black. }
\label{fig:OPF}
\end{figure}

We run our numerical experiments on a modified version of a 4-bus system in \cite{4bus} (case4gs in the Matpower library) which is illustrated in figure \ref{fig:OPF}.
The system has two conventional generators at Bus\,1 and Bus\,4, with active and reactive power limits $p_{Gi}^{min}=0,p_{Gi}^{max}=500$ and $q_{Gi}^{min}=-250,q_{Gi}^{max}=500$. Bus\,1 is the reference bus, while all other buses are PQ buses. We assume that the load at Bus\,2 is uncertain, with active power fluctuations $\omega$ uniformly distributed on $[-50,50]$. The reactive power fluctuations on Bus\,2 are proportional to the active power fluctuations, with $\gamma \approx 0.62$. We assume quadratic cost for Bus\,1 with $(c_{2,1},c_{1,1},c_{0,1})=(0.01,30,200)$ and a linear cost for Bus\,4 with $(c_{2,4},c_{1,4},c_{0,4})=(0,25,400)$.

\subsection{Numerical results}
We verify the approximation quality of the chance constraint approximation, and then assess the performance of the full CC-AC-OPF problem.

\subsubsection{Approximation of chance constraints}
We employ the two-step approach described in \ref{sec:two steps} to obtain the chance constraint approximations through the polynomials $h^\ast_0,\ldots,h^\ast_k$ given in \eqref{eq:acc-opnf}.
We investigate the accuracy of this approximation and how the accuracy improves by increasing the relaxation order $d$ and the addition of Stokes constraints. We show results for both the outer approximation \ref{sec:outer} and the (approximate) inner approximation \ref{sec:inner}.

To obtain outer and inner approximations we need to compute the probability of the projections of the sets $K_j$ defined in \eqref{eq:kj outer} and \eqref{eq:kj inner} respectively, by using the two-step method in Section~\ref{sec:two steps}. For the corresponding SDP relaxations, we choose the relaxation order of the first step to be $d=2$ or $3$ and for the second step to be $d+5=7$ or $8$.  For the first step, a lower degree polynomial is sufficient to approximate the level sets of $K_j$, whereas in the second step needs higher orders for better approximation and benefiting from Stokes constraints.

To assess how close we are to the true feasible set of the chance constraints, we created a large number of grid point to represent $B_\x$b using $100$ grid points for both active and reactive power for a total of $10'000$ grid points. For each grid point, we sampled 1'000 realizations of $\w$. For each $(\x,\w)$, we solved a standard power flow using Matpower. 
We then calculated the probability that a constraint holds for fixed \x\ by dividing the number of samples for which the power flow satisfies the constraints by the total number of samples for \w. 

\begin{figure}
\begin{center}
\includegraphics[width=0.8\columnwidth]{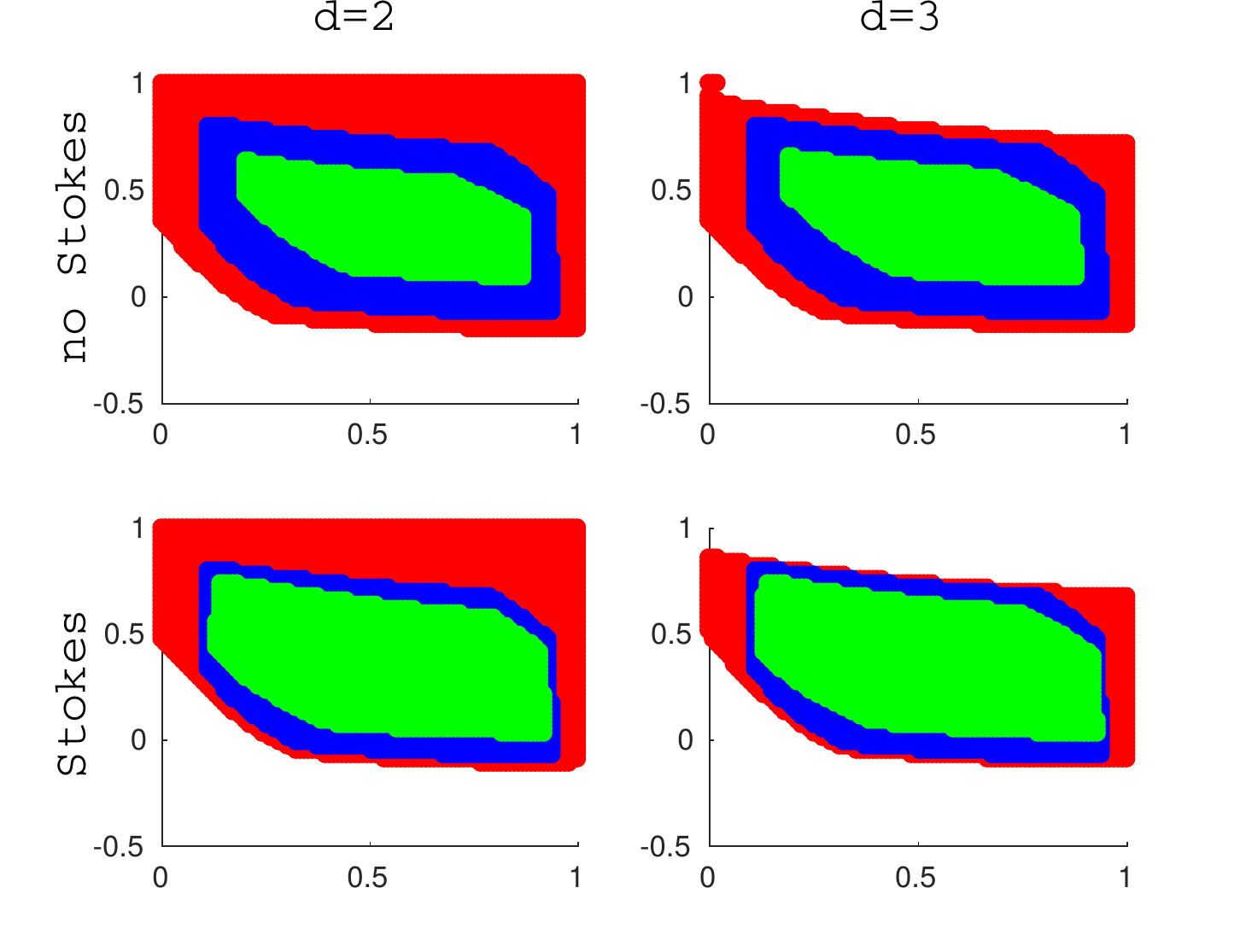}

\vspace{-15pt}
\caption{Comparison of the outer (red) and inner (green) approximation with the Monte Carlo simulation (blue) for $\varepsilon_1=0.01$ and $\varepsilon_2=0.1$.}
\label{fig:comparison eps01}
\end{center}
\end{figure}

Figure \ref{fig:comparison eps01} shows the feasible region for  $\varepsilon_1=0.01$ and $\varepsilon_2=0.1$. We show both the inner (green) and outer (red) approximation of the feasible region for relaxation orders $d=2,3$, and both with and without Stokes constraints. As a benchmark, we also show the feasible region computed through the Monte Carlo simulation (blue). The closer the approximated regions (green and red) are to the benchmark (blue), the better the approximation. We remark that both increasing the relaxation order and introducing Stokes constraints increase the quality of the solution. The improvement obtained by introducing Stokes constraints is very significant, while increasing the relaxation order only slightly increases the quality of the approximation. 


To further assess the quality of approximation, we report the ratios between the volume of the approximated feasibility regions and the volume computed through the Monte Carlo simulation in Table \ref{tab:ratio} for $\varepsilon_1=0.01$ and different values of $\varepsilon_2$. The addition of Stokes constraints clearly offers significant improvement. Interestingly, the quality of the outer approximation does not seem to depend on the choice of $\varepsilon_2$, while the accuracy of the inner approximation decreases with $\varepsilon_2$.

We observe that the outer approximation to the chance constraints is not very tight, and might lead to violation probabilities significantly above the acceptable levels. The extension proposed in this paper to allow for an (approximate) inner approximation provides a significant practical advantage over the previously existing methods in terms of returning safe approximations. It is also accurate enough to provide non-empty feasible sets, even at low relaxation orders.

\begin{table}
\begin{center}
\def\arraystretch{1.25}
\begin{tabular}{|c|cc|cc|cc|cc|}
\hline
                & \multicolumn{4}{c}{outer}&\multicolumn{4}{c|}{inner}\\
                & \multicolumn{2}{c}{$d=2$}&  \multicolumn{2}{c|}{$d=3$}& \multicolumn{2}{c}{$d=2$}&  \multicolumn{2}{c|}{$d=3$}\\
$\varepsilon_2$ &   --       &   Stokes      &   --       &   Stokes     &   --       &   Stokes      &   --       &   Stokes       \\
\hline
0.20            &   175\%   &   165\%   &  141\%    &  124\%     &    53\%    &   79\%   & 56\%   &    79\%  \\
0.15            &   179\%   &   168\%   &  143\%    &  126\%     &    49\%    &   75\%   & 53\%   &    76\%  \\
0.10            &   182\%   &   171\%   &  144\%    &  126\%     &    43\%    &   69\%   & 47\%   &    70\%  \\
0.05            &   185\%   &   173\%   &  144\%    &  125\%     &    30\%    &   56\%   & 36\%   &    58\%   \\
\hline
\end{tabular}
\vspace{6pt}
\caption{Ratio approximated vs. real volume for different values of $\varepsilon_2$.}
\label{tab:ratio}
\vspace{-24pt}
\end{center}
\end{table}

\subsection{Solving an instance of a CC-AC-OPF}
We assess the performance of the ACC-OPNF formulation in \eqref{eq:acc-opnf} by evaluating the cost of the optimal generation dispatch, the empirical constraint violation probability and by relating it to the deterministic AC-OPF.
For this experiment, we use the best inner approximation with relaxation order $d=3$ as well as the Stokes constraints to approximate the CC-AC-OPF \eqref{eq:abstract CCOPF}. We solve both the deterministic AC-OPF  and the approximation of the CC-AC-OPF  for different values of $\varepsilon_2$. We then compare the power injections, the cost and the maximal empirical violation probability of the individual chance constraints $\varepsilon_2^\ast$, which is computed through another Monte Carlo simulation at the obtained solution point using 1'000 samples of $\w$. 

Table \ref{tab:cc-ac-opf} summarizes the results. Column \emph{Det.} we show the results for the deterministic AC-OPF. The other columns are labeled by their acceptable violation probability for the individual constraint violation $\varepsilon_2$. The violation probability $\varepsilon_1=0.01$ for all experiments. 
The variables $p_{G0,4}$ and $q_{G0,4}$ are the independent variables $\x$ in our problem formulation, corresponding to the active and the reactive power of the generator at Bus $4$ in the test case.
The power injections at the slack bus generator $p_{G0,1}$ and $q_{G0,1}$ are among the dependent $\y_\x$ variables. Since these generators will adjust their values based on the realization of $\w$, we report their expected values in the table.
Further, we list the cost of the operating point and the maximum empirical violation probabilities $\epsilon_2^\ast$ among all individual constraints. We do not show results for the empirical violation probability of the joint chance constraint $\varepsilon^\ast_1$, as it was constantly $0\%$ for all optimal operating points. This is expected, since the engineering limits are typically more limiting than the power flow solvability conditions. 

As the violation probability $\epsilon_2$ decreases, more and more of the system load must be covered by the more expensive slack generator, resulting in a higher value for $p_{G0,1}$ and a higher expected cost. 
Considering the violation probabilities of the individual chance constraints we see that the optimal solution to the deterministic AC-OPF violates at least one of these constraints with a probability of almost $40\%$. For the approximations of the CC-AC-OPF the empirical violation probability $\varepsilon_2^\ast$ of the individual chance constraints is always below the requested probability $\varepsilon_2$, reflecting the fact that we indeed obtain a true inner approximation. While the empirical violation probability is quite close to the acceptable level for $\varepsilon_2=20\%$ and $\varepsilon_2=15\%$, respectively, the approximation is significantly more conservative for lower values of $\varepsilon_2$. For $\varepsilon_2=5\%$ no violations are observed.  

\begin{table}
\begin{center}
\def\arraystretch{1.5}
\begin{tabular}{|l|r|r|r|r|r|}
\hline
                 & Det.\hspace{0.5cm} & $\varepsilon_2= 20\%$&  $\varepsilon_2=15\%$ & $\varepsilon_2=10\%$& $\varepsilon_2=5\%$\\
\hline
$p_{G0,1}$       &    8.5   &    30.1    &    36.3    &   44.7    &   58.8 \\
$q_{G0,1}$       &  158.4   &   168.0    &   168.2    &  168.6    &  169.1 \\
$p_{G0,4}$       &  500.0   &   477.6    &   471.2    &  462.4    &  447.9 \\
$q_{G0,4}$       &  149.5   &   135.4    &   134.0    &  132.1    &  129.1 \\ 
\hline
cost             & 13\,357  &  13\,452   &   13\,481  &  13\,523  &  13\,596  \\
$\varepsilon_2^\ast$& 39.8\%   &  18.2\%    &   12.1\%   &  3.7\%    &    0.0\%  \\
\hline
\end{tabular}
\vspace{6pt}
\caption{Optimal values and solutions to \eqref{eq:abstract OPF} and \eqref{eq:acc-opnf} for $\varepsilon_1=0.01$ and different values of $\varepsilon_2$}
\vspace{-24pt}
\label{tab:cc-ac-opf}
\end{center}
\end{table}


\section{Conclusion}\label{sec:concl}

In this paper, we develop a new approach to handle chance constrained optimization problems in non-linear physical networks. 
The method is based on Semidefinite Programming (SDP) techniques to compute the volume of semi-algebraic sets, from which polynomial approximations of the chance constraints are obtained. To make existing results applicable in our practical setting, we (i) propose a set reformulation in order to enable inner approximations, and (ii) develop a two-step procedure to improve approximation quality at lower computational overhead. 

The method is applicable to any problem with polynomial equality and inequality constraints, and we demonstrate it numerically on the chance constrained AC Optimal Power Flow. In our experiments, the polynomial approximations were shown to provide sufficiently accurate representations of the feasible domain, and the resulting CC-AC-OPF was able to provide safe operating points with limited violation probability.

The method presented is a powerful and novel technique to handle chance constrained optimization for non-linear systems. Although, in its current form the method is applicable to small systems, it has the potential for multiple extensions and improvements. One promising future direction is to exploit the sparsity structure of networks to scale the method to larger instances.


\section*{Acknowledgment}
The first author is grateful to the Los Alamos National Laboratory for hosting him during summer 2017. His work is supported by ERC-Advanced Grant \#666981 TAMING.


\bibliographystyle{IEEEtran}
\bibliography{conestillmann}

 \appendix

\subsection{Proof}
\begin{proof} For the proof it will be handy to introduce some formulas. Define
 \[
 \begin{split}
 A:=&\exists \y\in Y, f(\x,\y,\w)= 0,\\
 B:=&\exists \y \in Y, (f(\x,\y,\w)= 0\land g_j(\x,\y,\w)\leq0,\\
 B':=&\forall \y \in Y, (f(\x,\y,\w)\neq 0\lor g_j(\x,\y,\w)>0),\\
 B'':=&\forall \y \in Y, (f(\x,\y,\w)\neq 0\lor g_j(\x,\y,\w)\geq0),\\
 C:=&\exists \y \in Y, (f(\x,\y,\w)= 0\land g_j(\x,\y,\w)\geq0).
 \end{split}
 \]
 Note that $\neg B = B'$ and $B'\Rightarrow B''$. Therefore $\x\in\kx$ is a stronger condition than 
 \[
 \begin{split}
 \x\in\kxpp&:= \{\x\in B_\x : \prob(A) \geq 1-\varepsilon_1,\\ 
 &\prob(B'') \geq 1-\varepsilon_2+ \varepsilon_1,\; j=1,\ldots,k\},
 \end{split}
 \]
 i.e., $\kx\subseteq\kxpp$. To see that $\kxpp\subseteq\kxp$, note that $B'' \Leftrightarrow (A\land B'') \lor (\neg A\land B'')$, $(A\land B'')\Rightarrow C$, and $(\neg A\land B'')\Leftrightarrow \neg A$. Hence, if $x\in \kx$, $1- \varepsilon_2 + \varepsilon_1 \leq \prob(B) \leq \prob(\neg A) + \prob(C)\leq \varepsilon_1 + \prob(C)$. Consequently, $ \prob(C)\geq 1-\varepsilon_2$, i.e., $x\in \kxp$.
\end{proof}

\end{document}